\documentclass{article}

\usepackage{amssymb, latexsym}

\usepackage{graphicx}

\usepackage{amsmath}

\newtheorem{theorem}{Theorem}[section]

\newtheorem{conjecture}[theorem]{Conjecture}

\newtheorem{corollary}[theorem]{Corollary}

\newtheorem{lemma}[theorem]{Lemma}

\newtheorem{proposition}[theorem]{Proposition}

\newtheorem{remark}[theorem]{Remark}

\begin{document}

\title{\Large\bf On the $L^p$-estimates for Beurling-Ahlfors  and Riesz transforms on Riemannian manifolds}

\author{Xiang-Dong Li\thanks{Research supported by NSFC No. 10971032, Key Laboratory RCSDS, CAS, No. 2008DP173182, and a
Hundred Talents Project of AMSS,
CAS.}\\
\\
{\it\small Academy of Mathematics and Systems Science, Chinese
Academy of Sciences}\\
{\it\small 55, Zhongguancun East Road, Beijing, 100190, P. R. China}\\
{\it\small E-mail: xdli@amt.ac.cn}
}

\maketitle

\begin{center}
\begin{minipage}{120mm}
\begin{center}{\bf Abstract}\end{center} In our previous papers \cite{Li2008, Li2011}, we
proved some martingale transform representation formulas for the
Riesz transforms and the Beurling-Ahlfors transforms on complete
Riemannian manifolds, and proved some explicit $L^p$-norm estimates
for these operators on complete Riemannian manifolds with suitable
curvature conditions. In this paper we correct a gap contained in
\cite{Li2008, Li2011} and prove that the $L^p$-norm of the
Riesz transforms $R_a(L)=\nabla(a-L)^{-1/2}$ can be explicitly bounded by $C(p^*-1)^{3/2}$ if $Ric+\nabla^2\phi\geq -a$ for $a\geq 0$, and the $L^p$-norm of the
Riesz transform $R_0(L)=\nabla(-L)^{-1/2}$ is bounded by $2(p^*-1)$ if $Ric+\nabla^2\phi=0$. We also prove that the $L^p$-norm estimates for the Beurling-Ahlfors transforms obtained in \cite{Li2011} remain valid. Moreover, we prove the time reversal martingale transform representation formulas for the Riesz transforms and the Beurling-Ahlfors transforms on complete  Riemannian manifolds.
\end{minipage}
\end{center}

\section{Introduction}

In our previous paper \cite{Li2008}, the author obtained a
martingale transform representation formula for the Riesz transforms
on complete Riemannian manifolds. More precisely, by the formula
(24) in Theorem 3.2 in \cite{Li2008}, the probabilistic
representation formula of the Riesz transform
$R_a(L)=\nabla(a-L)^{-1/2}$ acting on a nice function $f$ was given by
\begin{eqnarray*}
-{1\over 2}R_a(L)f(x)=\lim\limits_{y\rightarrow +\infty}E_y\left[
\left.  \int_0^\tau e^{a(s-\tau)}M_\tau M_s^{-1}dQ_af(X_s,
B_s)dB_s\right|X_\tau=x\right].\label{MRF1}
\end{eqnarray*}
Recently, R. Ba\~nuelos and F. Baudoin \cite{BB} pointed out that,
since $e^{-a\tau}M_\tau$ is not adapted to the filtration
$\mathcal{F}_t=\sigma(X_s, B_s, s\leq t)$, the above probabilistic
representation formula should be corrected as follows
\begin{eqnarray}
-{1\over 2}R_a(L)f(x)=\lim\limits_{y\rightarrow
+\infty}E_y\left[e^{-a\tau}M_\tau \left.  \int_0^\tau
e^{as}M_s^{-1}dQ_af(X_s, B_s)dB_s\right|X_\tau=x\right].\label{MRF2}
\end{eqnarray} Indeed, a careful check of the original proof of the formula $(24)$ in  Theorem
3.8 in \cite{Li2008} indicates that the correct probabilistic
representation formula of $R_a(L)f$ should be  given by
$(\ref{MRF2})$. See Section $2$ below. By the above observation, R. Ba\~nuelos and F.
Baudoin \cite{BB} pointed out that
 there is a gap in the proof of the $L^p$-norm
estimates of the Riesz transforms in
\cite{Li2008} and they proved a new martingale inequality which can be used to correct this gap. In this paper, we correct the above gap and prove that the $L^p$-norm of the Riesz transform $R_a(L)$ is bounded above by $C(p^*-1)^{3/2}$ if $Ric+\nabla^2\phi\geq -a$ for $a\geq 0$, and the $L^p$-norm of the Riesz transform $R_0(L)$ is bounded by $2(p^*-1)$ if $Ric+\nabla^2\phi=0$. See Theorem \ref{EST1} below. We also correct the gap contained in \cite{Li2011} (due to the same reason as above) and prove that the main results on the $L^p$-norm estimates of the Beurling-Ahlfors transforms obtained in \cite{Li2011} remain valid. See Theorem \ref{BA} and Remark \ref{rrr} below. Moreover, we prove the time reversal martingale transform representation formulas for the Riesz transforms and the Beurling-Ahlfors transforms on complete  Riemannian manifolds.

\section{Riesz transforms on functions}

Let $(M, g)$ be a complete Riemannian manifold, $\nabla$ the
gradient operator on $M$, $\Delta$ the Laplace-Beltrami operator on
$M$. Let $\phi\in C^2(M)$, and
$d\mu=e^{-\phi}dv$, where $dv$ is the standard Riemannian volume
measure on $M$.
Let $L_0^2(M, \mu)=L^2(M, \mu)$ if $\mu(M)=\infty$, and $L^2_0(M, \mu)=\{f\in  L^2(M, \mu):  \int_M f d\mu=0\}$ if $\mu(M)<\infty$.

Let $L=\Delta-\nabla\phi\cdot\nabla$. Let $d$ be the exterior differential operator, $d^*_\phi$ be its $L^2$-adjoint with respect to the weighted volume measure  $d\mu=e^{-\phi}dv$. Let $\square_\phi=dd_\phi^*+d_\phi^*d$ be the Witten-Laplacian acting on forms over $(M, g)$ with respect to the weighted volume measure $d\mu=e^{-\phi}dv$.

Let $B_t$ be one dimensional Brownian motion on $\mathbb{R}$
starting from $B_0=y>0$ and with infinitesimal generator ${1\over
2}{d^2\over dy^2}$. Let
$$
\tau=\inf\limits\{t>0: B_t=0\}.$$

Let $X_t$ be the $L$-diffusion process on $M$.   Let $Ric$ be the Ricci curvature
on $(M, g)$, $\nabla^2\phi$ be the Hessian of the potential function
$\phi$. Let $M_t\in {\rm
End}(T_{X_0}M, T_{X_t}M)$ is the unique solution to the covraiant SDE along the trajectory of $(X_t)$:
\begin{eqnarray*}
{\nabla \over \partial t}M_t=-(Ric+\nabla^2\phi)(X_t)M_t, \ \ \ \
M_0={\rm Id}_{T_{X_0}M}.
\end{eqnarray*}
In particular, in the case where $Ric+\nabla^2\phi=-a$, we have
$$M_t=e^{at}U_{t}, \ \ \ \ \ \forall t\geq 0,$$
where $U_t: T_{X_0}M\rightarrow T_{X_t}M$ denotes the stochastic parallel
transport along $X_t$.

The following result is the correct reformulation of Lemma 3.7 in \cite{Li2008}.

\begin{lemma}\label{lem1} For all $\eta\in C_0^\infty(M, \Lambda^1 T^*M)$, and $\eta_a(x, y)=e^{-y\sqrt{a+\square_\phi}}\eta(x)$, we have
\begin{eqnarray}
\eta(X_\tau)=e^{a\tau}M_{\tau, k}^{*, -1}\eta_a(X_0, B_0)+e^{a\tau}M_{\tau, k}^{*}\int_0^\tau e^{-as}M_{s}^*\left(\nabla, {\partial\over \partial y}\right)\eta_a(X_s, B_s)\cdot(U_sdW_s, dB_s).
\label{w1}
\end{eqnarray}
\end{lemma}
{\it Proof}. By It\^o's calculus, we have (see p.266 line 16 in \cite{Li2008})
\begin{eqnarray*}
{\nabla\over \partial t}\left(e^{-at}M_{t}^*\eta_a(X_t, B_t)\right)=e^{-at}M_{t}^*\left(\nabla, {\partial\over \partial y}\right)\eta_a(X_t, B_t)\cdot(U_tdW_t, dB_t).
\end{eqnarray*}
Integrating from $t=0$ to $t=\tau$ , we complete the proof of Lemma \ref{lem1}. \hfill $\square$

\medskip

The following result is the correct reformulation of Theorem 3.8 in \cite{Li2008}.

\begin{theorem}\label{theo1} Let $\omega\in C_0^\infty(M< \Lambda^1T^*M)$, and $\omega_a(x, y)=e^{-y\sqrt{a+\square_\phi}}\omega(x)$. Then
\begin{eqnarray}
{1\over 2}\omega(x)=\lim\limits_{y\rightarrow \infty}E_y\left[e^{-a\tau}M_{\tau}\left.
\int_0^\tau e^{as}M_{s}^{-1}{\partial \over \partial y}\omega_a(X_s, B_s)dB_s\right|X_\tau=x\right].\label{w2}
\end{eqnarray}
\end{theorem}
{\it Proof}. The proof is indeed a small modification of the original proof of Theorem 3.8 given in \cite{Li2008}. For the completeness of the paper, we produce the details here. Let $Z_t=(X_t, B_t)$, $\eta\in C_0^\infty(\Lambda^kT^*M)$. By $(\ref{w1})$ in Lemma \ref{lem1}, we have
\begin{eqnarray*}
\eta(X_\tau)=e^{a\tau}M_\tau^{*, -1}\eta_a(Z_0)+e^{a\tau}M_{\tau}^{*, -1}\int_0^\tau e^{-as}M_s^*\left(\nabla, {\partial\over \partial y}\right)\eta_a(Z_r)\cdot(U_sdW_s, dB_s).
\end{eqnarray*}
Hence
\begin{eqnarray*}
& &\int_M\left\langle E_y\left[e^{-a\tau}M_{\tau}\left.\int_0^\tau e^{as}M_{s}^{-1}{\partial \over \partial y}\omega_a(X_s, B_s)dB_s\right|X_\tau=x\right], \eta(x)\right\rangle d\mu(x)\\
& &\ \ \ = E_y\left[e^{-a\tau}M_{\tau}\left.\int_0^\tau e^{as}M_{s}^{-1}{\partial \over \partial y}\omega_a(X_s, B_s)dB_s, \eta(X_\tau)\right\rangle\right]\\
& &\ \ \ =I_1+I_2,
\end{eqnarray*}
where
\begin{eqnarray*}
I_1&=&E_y\left[\left\langle e^{-a\tau}M_{\tau}\int_0^\tau e^{as}M_{s}^{-1}{\partial \over \partial y}\omega_a(X_s, B_s)dB_s, e^{a\tau}M_{\tau}^{*, -1}\eta_a(X_0, B_0)\right\rangle\right],\\
I_2&=&E_y\left[\left\langle e^{-a\tau}M_{\tau}\int_0^\tau e^{as}M_{s}^{-1}{\partial \over \partial y}\omega_a(X_s, B_s)dB_s,\right.\right. \\
& & \ \ \ \ \ \ \ \ \ \ \ \ \ \ \   \left.\left. e^{a\tau}M_{\tau}^{*, -1}\int_0^\tau e^{-as}M_{s}^{*}(\nabla, \partial_y)\eta_a(X_s, B_s)\cdot (U_s dW_s, dB_s)\right\rangle\right].
\end{eqnarray*}
Using the martingale property of the It\^o integral, we have
\begin{eqnarray*}
I_1&=&E_y\left[\left\langle \int_0^\tau e^{as}M_{s}^{-1}{\partial \over \partial y}\omega_a(X_s, B_s)dB_s, \eta_a(X_0, B_0)\right\rangle\right]\\
&=&E_y\left[\left\langle E\left[ \left.\int_0^\tau e^{as}M_{s}^{-1}{\partial \over \partial y}\omega_a(X_s, B_s)dB_s\right|(X_0, B_0)\right], \eta_a(X_0, B_0)\right\rangle\right]\\
&=&0.
\end{eqnarray*}
On the other hand, using the $L^2$-isometry of the It\^o integral, we have
\begin{eqnarray*}
I_2&=&E_y\left[\left\langle \int_0^\tau e^{as}M_{s}^{-1}{\partial \over \partial y}\omega_a(X_s, B_s)dB_s, \int_0^\tau e^{-as}M_{s}^{*}(\nabla, \partial_y)\eta_a(X_s, B_s)\cdot (U_s dW_s, dB_s)\right\rangle\right]\\
&=&E_y\left[\int_0^\tau \left\langle  e^{as}M_{s}^{-1}{\partial \over \partial y}\omega_a(X_s, B_s), e^{-as}M_{s}^{*}{\partial\over \partial y}\eta_a(X_s, B_s)\right\rangle ds\right]\\
&=&E_y\left[\int_0^\tau \left\langle  {\partial \over \partial y}\omega_a(X_s, B_s), {\partial\over \partial y}\eta_a(X_s, B_s)\right\rangle ds\right].
\end{eqnarray*}
The Green function of the background radiation process is given by $2(y\wedge z)$. Thus
\begin{eqnarray*}
& &E_y\left[\int_0^\tau \left\langle  {\partial \over \partial y}\omega_a(X_s, B_s), {\partial\over \partial y}\eta_a(X_s, B_s)\right\rangle ds\right]\\
& &\ \ \ \ \ \ \ \ \ =2\int_M\int_0^\infty (y\wedge z)\left\langle {\partial\over \partial z}\omega_a(x, z),  {\partial\over \partial z}\eta_a(x, z)\right\rangle dzd\mu(x).
\end{eqnarray*}
By spectral decomposition, we have the Littelwood-Paley identity
\begin{eqnarray*}
\lim\limits_{y\rightarrow \infty}\int_M\int_0^\infty (y\wedge z)\left\langle {\partial\over \partial z}\omega_a(x, z),  {\partial\over \partial z}\eta_a(x, z)\right\rangle dzd\mu(x)=\int_M \langle \omega(x), \eta(x)\rangle d\mu(x).
\end{eqnarray*}
Thus
\begin{eqnarray*}
\langle \omega, \eta\rangle_{L^2(\mu)}=2\lim\limits_{y\rightarrow \infty} \int_M \left\langle E_y\left[ e^{-a\tau}M_{\tau}\left.
\int_0^\tau e^{as}M_{s}^{-1}{\partial\over \partial y}\omega_a(X_s, B_s)dB_s\right|X_\tau=x\right], \eta(x)\right\rangle d\mu(x).
\end{eqnarray*}
This completes the proof of Theorem \ref{theo1}. \hfill $\square$

\medskip

The following martingale transform representation formula of the Riesz transforms on complete Riemannian manifolds, which is the extension of
 the Gundy-Varopoulos representation formula
of the Riesz transforms on Euclidean space \cite{GV}, is the correct reformulation of the one that we obtained in Theorem 3.2 in \cite{Li2008}.

\begin{theorem}\label{Th-2} Let $R_a(L)=\nabla(a-L)^{-1/2}$. Then, for all $f\in C_0^\infty(M)$, we have
\begin{eqnarray}
R_a(L)f(x)=-2\lim\limits_{y\rightarrow +\infty}E_y\left[
\left. e^{-a\tau}M_\tau \int_0^\tau e^{as}M_s^{-1}dQ_af(X_s,
_s)dB_s\right|X_\tau=x\right].\label{R1}
\end{eqnarray} In particular, in the case where $Ric+\nabla^2\phi=-a$, we have
\begin{eqnarray}
R_a(L)f(x)=-2\lim\limits_{y\rightarrow
+\infty}E_y\left[\left. U_\tau \int_0^\tau U_s^{-1}dQ_af(X_s,
B_s)(U_sdW_s, dB_s)\right|X_\tau=x\right].\label{R2}
\end{eqnarray}
\end{theorem}
{\it Proof}. Applying Theorem \ref{theo1} to $\omega=d(a-L)^{-1/2}f$, the proof of Theorem \ref{Th-2} is as the same as the one of Theorem 3.2 given in \cite{Li2008}. \hfill $\square$

We now state the $L^p$-norm estimates of the Riesz transforms on complete Riemannian manifolds.
Throughout this paper, for any $p\in (1, \infty)$, let
$$
p^*=\max\limits\left\{p, \ {p\over p-1}\right\}.
$$

The following result is a correction of Theorem 1.4 in \cite{Li2008}.

\begin{theorem}\label{EST1} Let $M$ be a complete Riemannian manifold, and $\phi\in
C^2(M)$. Then \\
(i)\ \ for all $f\in C_0^\infty(M)$,
\begin{eqnarray}
\|\nabla(a-L)^{-1/2}f\|_2\leq \|f\|_2,
\end{eqnarray}
(ii)\ \ if $Ric+\nabla^2\phi\equiv 0$, then for all $p\in (1, \infty)$,
\begin{eqnarray} \|\nabla(-L)^{-1/2}f\|_p\leq 2(p^*-1)\|f\|_p, \ \ \
\forall f\in C_0^\infty(M), \ \label{a1}
\end{eqnarray}
if $Ric+\nabla^2\phi\equiv -a$, where $a>0$ is a constant, then for all $p\in (1, \infty)$,
\begin{eqnarray}
\|\nabla(a-L)^{-1/2}f\|_p\leq 2(p^*-1)(1+4\|T_1\|_p)\|f\|_p, \ \ \
\forall f\in C_0^\infty(M), \ \label{a2}
\end{eqnarray}
where $T_1$ is the first exiting time of the standard
$3$-dimensional Brownian motion from the unit ball $B(0, 1)=\{x\in
\mathbb{R}^3: \|x\|=1\}$.\\
(iii) if $Ric+\nabla^2\phi\geq -a$, where $a\geq 0$ is a constant,
then there is a numerical constant $C>0$ such that for all $p>1$,
\begin{eqnarray}
\|\nabla(a-L)^{-1/2}f\|_p\leq C(p^*-1)^{3/2}\|f\|_p, \ \ \ \forall
f\in C_0^\infty(M). \ \label{a3}
\end{eqnarray}
\end{theorem}
{\it Proof}. The case $(i)$ for $p=2$ is well known, cf. \cite{Li2008, Li2010}. By \cite{Li2008}, for any fixed $x\in M$, there exists a bounded operator $A(x)\in {\rm End}(T_{x}M)$ such that that $d\omega(x)=A\nabla \omega(x)$ and $\|A(x)\|_{\rm op}\leq 1$. In the case $Ric+\nabla^2\phi=-a$, we have
\begin{eqnarray*}
\nabla(a-L)^{-1/2}f(x)=-2\lim\limits_{y\rightarrow
+\infty}E_y\left[\left. U_\tau \int_0^\tau U_s^{-1}A\nabla Q_af(X_s,
B_s)dB_s\right|X_\tau=x\right].
\end{eqnarray*}
The stochastic integral in the above formula is a subordination of martingale transforms. By Burkholder's sharp $L^p$-inequality for
martingale transforms \cite{Bk} we obtain
\begin{eqnarray*}
\|\nabla(a-L)^{-1/2}f\|_p\leq 2(p^*-1)\sup\limits_{s\in [0,
\tau]}\|A(X_s)\|_{\rm op}\left\|\int_0^\tau (\nabla,
\partial_y)Q_a(f)(X_s, B_s)\cdot (U_sdW_s, dB_s)\right\|_p,
\end{eqnarray*}
where $\|A(X_s)\|_{\rm op}$ denotes the operator norm of $A(X_s)$ on
$T_{X_s}M$. Note that  $$\sup\limits_{s\in [0, \tau]}\|A(X_s)\|_{\rm op}\leq 1.$$
This yields
\begin{eqnarray*}
\|\nabla(a-L)^{-1/2}f\|_p\leq 2(p^*-1)\left\|\int_0^\tau (\nabla,
\partial_y)Q_a(f)(X_s, B_s)\cdot (U_sdW_s, dB_s)\right\|_p,
\end{eqnarray*}
In \cite{Li2008}, we have proved that, for all $1<p<\infty$, it holds
\begin{eqnarray*}
\left\|\int_0^\tau (\nabla,
\partial_y)Q_a(f)(X_s, B_s)\cdot (U_sdW_s, dB_s)\right\|_p\leq \left(1+4\|T_1\|_p1_{a>0}\right)\|f\|_p.
\end{eqnarray*}
Combining this with the previous inequality, we obtain
\begin{eqnarray*}
\|\nabla(a-L)^{-1/2}f\|_p\leq 2(p^*-1)(1+4\|T_1\|_p1_{a>0})\|f\|_p.
\end{eqnarray*}
This proves the case of (ii).
\medskip

In general case $Ric+\nabla^2\phi\geq -a$, we have
\begin{eqnarray*}
\nabla(a-L)^{-1/2}f(x)=2\lim\limits_{y\rightarrow
+\infty}E_y\left[\left. e^{-a\tau} M_\tau \int_0^\tau e^{as} M_s^{-1}A\nabla Q_af(X_s,
B_s)dB_s\right|X_\tau=x\right].
\end{eqnarray*}
By the $L^p$-contractivity of conditional expectation, see \cite{Li2008}, we have
\begin{eqnarray*}
\|\nabla(a-L)^{-1/2}f\|_p\leq 2\lim\inf\limits_{y\rightarrow \infty}\left\|e^{-a\tau}M_\tau \int_0^\tau e^{as} M_s^{-1}A\nabla Q_af(X_s,
B_s)dB_s\right\|_p.
\end{eqnarray*}
Let
$$
J_y=\left\{\int_0^\tau |\nabla Q_a f(X_s,
B_s)|^2ds\right\}^{1/2}.$$
By Theorem 2.6 due to Ba\~nuelos and Baudoin in \cite{BB}, under the condition $Ric+\nabla^2\phi\geq -a$, we can prove that
\begin{eqnarray*}
\left\|e^{-a\tau} M_\tau\int_0^\tau e^{as}M_s^{-1}dQ_af(X_s,
B_s)dB_s\right\|_p\leq 3\sqrt{p(2p-1)}\|J_y\|_p.
\end{eqnarray*}
By Proposition 6.2 in our previous paper \cite{Li2010}, for all $p\in (1, \infty)$,  we proved that
$$
\|J_y\|_p\leq B_p\|f\|_p,$$
where for all $p\in (1, 2)$, $
B_p=(2p)^{1/2}(p-1)^{-3/2}$, $B_2=1$, and for all $p\in (2, \infty)$, $B_p={p\over \sqrt{2(p-2)}}$. From the above  estimates, for all $p\in (1, 2)$, we can obtain
\begin{eqnarray*}
\|\nabla(a-L)^{-1/2}f\|_p&\leq &6\sqrt{2}p^{3/2}(2p-1)^{1/2}(p-1)^{-3/2}\|f\|_p\\
&\leq& 12\sqrt{6}(p-1)^{-3/2}\|f\|_p,
\end{eqnarray*}
and for $p>2$,
\begin{eqnarray*}
\|\nabla(a-L)^{-1/2}f\|_p
&\leq& 3\sqrt{2}p^{3/2}(2p-1)^{1/2}(p-2)^{-1/2}\|f\|_p\\
&\leq& 6(p-1)^{3/2}(1+O(1/p))\|f\|_p.
\end{eqnarray*}
The proof of Theorem \ref{EST1} is completed. \hfill $\square$

\begin{remark}{\rm The above proof corrects a gap in the proof of Theorem 1.4 given in \cite{Li2008} (p.270 line 9 to line 12 in \cite{Li2008}), where we used the Burkholder sharp $L^p$-inequality for martingale transforms. As $e^{-a\tau}M_{\tau}$ is not adapted with respect to the filtration $\mathcal{F}_s=\sigma(X_u, B_u, u\in [0, s])$, $s<\tau$, the proof given in \cite{Li2008} is valid only in the case $e^{-a\tau}M_\tau$ is independent of $(X_s: s\in [0, \tau])$, which only happens if $Ric+\nabla^2\phi\equiv -a$ for some constant $a\geq 0$.}.
\end{remark}

The following result is the correction of  Corollary 1.5 in \cite{Li2008}.

\begin{corollary}\label{EST2} Let $M$ be a complete Riemannian manifold with
non-negative Ricci curvature. Then there exists a numerical constant $C>0$ such that for all $p>1$,
\begin{eqnarray*}
\|\nabla(-\Delta)^{-1/2}f\|_p\leq C(p^*-1)^{3/2}\|f\|_p. \label{a4}
\end{eqnarray*}
In particular, if $Ric=0$, i.e., if $M$ is a Ricci flat Riemannian
manifold, then for all $1<p<\infty$,
\begin{eqnarray*}
\|\nabla(-\Delta)^{-1/2}f\|_p\leq 2(p^*-1)\|f\|_p. \label{a5}
\end{eqnarray*}
\end{corollary}

In view of Theorem \ref{EST1} and Corollary \ref{EST2}, we need to reformulate Conjecture 1.7 in \cite{Li2008} as follows.

\begin{conjecture} Let $M$ be a complete Riemannian manifold, $\phi\in C^2(M)$. Suppose
that $Ric(L) = Ric +\nabla^2\phi=0$. Then there exists a constant $c > 0$ such that for all
$p > 1$, we have
\begin{eqnarray*}
c(p^*-1)(1 + o(1))\leq \|\nabla(-L)^{-1/2}\|_{p, p}\leq 2(p^*-1).
\end{eqnarray*}
In particular, on any complete Riemannian manifold M with flat Ricci
curvature, for all $p > 1$, we have
\begin{eqnarray*}
c(p^*-1)(1 + o(1))\leq \|\nabla(-\Delta)^{-1/2}\|_{p, p}\leq 2(p^*-1).
\end{eqnarray*}

\end{conjecture}

\begin{remark}{\rm Using the Bellman function technique,  Carbonaro and Dragi\v{c}evi\`{c} \cite{CD2012} proved that if $Ric+\nabla^2\phi\geq -a$, then for all $p\in (1, \infty)$,
\begin{eqnarray*}
\|\nabla(a-L)^{-1/2}f\|_p\leq 12(p^*-1)\|f\|_p, \ \ \ \forall f\in C_0^\infty(M).
\end{eqnarray*}
It would be nice if one can find a probabilistic proof of this result.
}
\end{remark}

\section{Riesz transforms on Gaussian spaces}

In this section, we give the proof of Corollary 1.6 in \cite{Li2008}. Let $G$ be a compact Lie group endowed with a bi-invariant Riemannian metric, $\mathcal{G}$ its Lie algebra, and $n={\rm dim}G$. Let $X_1, \ldots, X_n$ be an orthonormal basis of $\mathcal{G}$, and $\Delta_G=\sum\limits_{i=1}^n X_i^2$ the Laplace-Beltrami operator on $G$. In \cite{Ar}, Arcozzi proved that, the $L^p$-norm of the Riesz transform $R^G:=\sum\limits_{i=1}^n R_{X_i}X_i$ on $G$ satisfies
$\|R^G\|_{p}\leq 2(p^*-1)$ for all $p\in (1, \infty)$, where $R_{X_i}=X_i(-\Delta_G)^{-1/2}$ is the Riesz transform on $G$ in the direction $X_i$. As the unit sphere $S^{n-1}$ can be identified as $S^{n-1}=SO(n)/SO(n-1)$, where $SO(n)$ is the rotation group of $\mathbb{R}^n$, Arcozzi proved that the $L^p$-norm of the Riesz transform $R^{S^{n-1}}=\nabla^{S^{n-1}}(-\Delta_{S^{n-1}})^{-1/2}$ on $S^{n-1}$ satisfies $\|R^{S^{n-1}}\|_p\leq 2(p^*-1)$ for all $p\in (1, \infty)$. Let $S^{n-1}(\sqrt{n})$ be the $(n-1)$-dimensional sphere of radius $\sqrt{n}$. Then the $L^p$-norm of the Riesz transform $R^{S^{n-1}(\sqrt{n})}$ satisfies  $\|R^{S^{n-1}(\sqrt{n})}\|_p\leq 2(p^*-1)$.  By the Poincar\'e limit, as $n\rightarrow \infty$, $S^{n-1}(\sqrt{n})$ endowed with the normalized volume measure converges in a proper way to the infinite dimensional Wiener space $\mathbb{R}^\mathbb{N}$ endowed with the Wiener measure, and the Laplace-Beltrami operator on $S^{n-1}(\sqrt{n})$ converges to the Orisntein-Uhlenbeck operator on $\mathbb{R}^\mathbb{N}$. From this, Arcozzi derived that the Riesz transform associated with the Ornstein-Uhlenbeck operator $L=\Delta-x\cdot \nabla$ on the Wiener space satisfies $\|\nabla(-L)^{-1/2}\|_p\leq 2(p^*-1)$ for all $p\in (1, \infty)$.

In general, let $A\in M(n, \mathbb{R})$ be a positive definite symmetric matrix on $\mathbb{R}^n$, and let $\langle x, y\rangle_A=\langle x, Ay\rangle$, $\forall x, y\in \mathbb{R}^n$. Then $\sqrt{A}: (\mathbb{R}^n, \langle\cdot,\cdot\rangle)\rightarrow (\mathbb{R}^n, \langle\cdot, \cdot\rangle_A)$ is an isometry. Let $SO(n, A)$ be the rotation group on $(\mathbb{R}^n, \langle\cdot,\cdot\rangle_A)$, and $S^{n-1}_A$ be the $(n-1)$-dimensional sphere in $(\mathbb{R}^n, \langle\cdot,\cdot\rangle_A)$. Then $S^{n-1}_{A}={SO(n, A)/SO(n-1, A)}$. By the same argument as used by Arcozzi \cite{Ar}, we can prove that the $L^p$-norm of the Riesz transform
on $SO(n, A)$ satisfies $\|R^{SO(n, A)}\|_p\leq 2(p^*-1)$, and the $L^p$-norm of the Riesz transform
on $S^{n-1}_{A}$ satisfies $\|R^{S^{n-1}_{A}}\|_p\leq 2(p^*-1)$. Similarly, we have $\|R^{S^{n-1}_{A}(\sqrt{n})}\|_p\leq 2(p^*-1)$. Thus, we have proved the following

\begin{theorem} \label{Th-3}  Let $A\in M(n, \mathbb{R})$ be a positive definite symmetric matrix  on $\mathbb{R}^n$, and let
$$L_A=\Delta-Ax\cdot \nabla$$ be the Ornstein-Uhlenbeck operator on the Gaussian space $(\mathbb{R}^n, \mu_A)$, where
$$d\mu_A(x)={1\over (2\pi \ {\rm det} A)^{n/2}}e^{-\langle x, Ax\rangle}dx.$$ Then, for all $1<p<\infty$, the $L^p$-norm of the Riesz transform $R=\nabla(-L_A)^{-1/2}$ on $(\mathbb{R}^n, \mu_A)$ satisfies
$$\|\nabla(-L_A)^{-1/2}\|_p\leq 2(p^*-1).$$
\end{theorem}

Using the Poincar\'e limit, we can derive the following result from Theorem \ref{Th-3}.

\begin{theorem}\label{Th-4}  (i.e., Corollary 1.6 in \cite{Li2008}) Let $(W, H, \mu_A)$ be an abstract Wiener space, where $W$ is a real separable Banach space, $H$ is a real separable Hilbert space which is densely embedded in $W$, $A\in \mathcal{L}(H)$ be a self-adjoint positive definite operator with finite Hilbert-Schmidt norm, and $\mu$ the Gaussian measure on $W$ with mean zero and with covariance $A$. Let
$$L_A=\Delta-Ax\cdot \nabla$$ be the generalized Ornstein-Uhlenbeck operator on $(W, H, \mu_A)$. Then, for all $1<p<\infty$, the $L^p$-norm of the Riesz transform $R=\nabla(-L_A)^{-1/2}$ on $(W, H, \mu_A)$ satisfies
$$\|\nabla(-L_A)^{-1/2}\|_p\leq 2(p^*-1).$$
\end{theorem}

\section{Beurling-Ahlfors transforms}

Throughout this section, let $M$ be a complete and stochastically complete Riemannian manifold, $n={\rm dim} M$. Let $X_t$ be Brownian motion on $M$, $W_k$ the $k$-th Weitzenb\"ock
curvature operator.
Let $A_i\in {\rm End}(\Lambda^kT^*M)$, $i=1, 2$, be the bounded endomorphism which, in a local normal coordinate $(e_1, \ldots, e_n)$ at any fixed point $x$, is defined by
$$
A_1=(a_ia_j^*)_{n\times n}, \ \ \ \  A_2=(a_i^*a_j)_{n\times n},$$
where $\ a_i={\rm int}_{e_i}$ is the inner multiplication by $e_i$, and $a_j^*=e_j^*\wedge$ is the exterior multiplication by $e_j$, $i, j=1, \ldots, n$. For details, see \cite{Li2011}.

Let
$M_t\in {\rm End}(\Lambda^kT_{X_0}^*M, \Lambda^kT_{X_t}^*M)$ be defined by
$${\nabla M_t\over \partial t}=-W_k(X_t)M_t, \ \ M_0={\rm Id}_{\Lambda^kT^*_{X_0}M}.$$
For any fixed $T>0$, the backward heat semigroup generated by the Hodge Laplacian
$\square$ on $k$-forms is defined by
$$
\omega(x, T-s)=e^{-(T-s)\square}\omega(x), \ \ \ \forall x\in M, s\in [0, T], \ \omega\in C_0^\infty(\Lambda^kT^*M).$$
Recall that, the Weitzenb\"ock formula reads as follows
$$\square=-{\rm Tr}\nabla^2+W_k.$$

We now state the martingale transform representation formula for the Beurling-Ahlfors
transforms on $k$-forms over complete Riemannian manifolds .

\begin{theorem}\label{RPF1} Let $M$ be a complete and stochastically complete Riemannian manifold. Suppose that $W_k\geq -a$, where $a\geq 0$ is a constant. Then, for all $\omega, \eta\in C_0^\infty(\Lambda^kT^*M)$, we have
\begin{eqnarray*}
\langle dd^*(a+\square)^{-1}\omega, \eta\rangle=2\lim\limits_{T\rightarrow \infty}\int_M \langle S_{A_2}^T\omega, \eta\rangle dx,
\end{eqnarray*}
\begin{eqnarray*}
\langle d^*d(a+\square)^{-1}\omega, \eta\rangle=2\lim\limits_{T\rightarrow \infty}\int_M \langle S_{A_1}^T\omega, \eta\rangle dx,
\end{eqnarray*}
where, for a.s. $x\in M$,
\begin{eqnarray*}
S_{A_i}^T\omega(x)=E\left[M_T e^{-aT}\left.\int_0^T e^{at}M_t^{-1}A_i\nabla\omega_a(X_t, T-t)dX_t\right|X_T=x\right], \ \ \ i=1, 2.
\end{eqnarray*}
In particular, the Beurling-Ahlfors transform
$$
S_B\omega:=(d^*d-dd^*)(a+\square)^{-1}\omega
$$
has the following martingale transform representation: for a.s. $x\in M$,
\begin{eqnarray*}
S_B\omega(x)=2\lim\limits_{T\rightarrow \infty}E\left[M_T e^{-aT}\left.\int_0^T e^{at}M_t^{-1}B\nabla\omega_a(X_t, T-t)dX_t\right|X_T=x\right],
\end{eqnarray*}
where
$$B=A_1-A_2.$$
\end{theorem}

\begin{remark} \label{rem1} {\rm The martingale transform representation formulas  in Theorem \ref{RPF1} are the correct reformulation of the formulas that we obtained in Theorem 3.4 in \cite{Li2011}, where the martingale transform representation formulas of $S_{A_i}$ and $S_B$ were given in the following  way
\begin{eqnarray*}
S_{A_i}^T\omega(x)=E\left[\left.\int_0^T e^{a(t-T)}M_TM_t^{-1}A_i\nabla\omega_a(X_t, T-t)dX_t\right|X_T=x\right], \ \ \ i=1, 2,
\end{eqnarray*}
and
\begin{eqnarray*}
S_B\omega(x)=2\lim\limits_{T\rightarrow \infty}E\left[\left.\int_0^T e^{a(t-T)}M_TM_t^{-1}B\nabla\omega_a(X_t, T-t)dX_t\right|X_T=x\right].
\end{eqnarray*}
The same correction should also be made for Theorem 3.5 in \cite{Li2011}, where $a=0$. The reason is that, as pointed out by Ba\~nuelos and Baudoin in \cite{BB}, $M_T$ is not adapted with respect to the filtration  $\mathcal{F}_t=\sigma(X_s: s\in [0, t])$, $t<T$. Moreover, in the proof of Theorem 1.2 in \cite{Li2011} (p.135, line 7 to line 8), we used the Burkholder-Davis-Gundy inequality to derive that
\begin{eqnarray*}
\|S_{A_i}^T\omega\|_p\leq C_p\sup\limits_{0\leq t\leq T}\|e^{a(t-T)}M_TM_t^{-1}A_i\|_{\rm op}\left\|\left\{\int_0^T |\overline{\nabla}\omega_a(X_t, T-t)|^2dt\right\}^{1/2}\right\|_p,
\end{eqnarray*}
where $\|\cdot\|_{\rm op}$ denotes the operator norm, and $C_p$ is a constant.
However, except that $M_T$ is independent of the $(X_t: t\in [0, T])$, one cannot use the Burkholder-Davis-Gundy inequality in above way,  due to the fact that $M_T$ is not adapted with respect to the filtration $\mathcal{F}_t=\sigma(X_s: s\in [0, t])$, $t<T$.
}
\end{remark}

\noindent{\it Proof of Theorem \ref{RPF1}}. By Remark \ref{rem1}, we need only to correct the martingale transform representation formulas appeared in Theorem 3.4 and Theorem 3.5 in \cite{Li2011} in the right way stated in Theorem \ref{RPF1}. Thus, the original proof given in \cite{Li2011} for these  formulas remain valid after a small modification. To save the length of the paper, we omit it here. \hfill $\square$

\medskip

\begin{proposition} \label{BA2} For all constant $a\geq 0$ and $\omega\in C_0^\infty(\Lambda^kT^*M)$, we have
\begin{eqnarray}
\|dd^*(a+\square)^{-1}\omega\|_2^2+\|d^*d(a+\square)^{-1}\omega\|_2^2=\|\square(a+\square)^{-1}\omega\|_2^2,\label{ID1}
\end{eqnarray}
 Moreover,
\begin{eqnarray*}
\|dd^*(a+\square)^{-1}\omega\|_2&\leq& \|\omega\|_2,\\
\|d^*d(a+\square)^{-1}\omega\|_2&\leq& \|\omega\|_2,
\end{eqnarray*}
and
\begin{eqnarray*}
\|(d^*d-dd^*)(a+\square)^{-1}\omega\|_2\leq 2\|\omega\|_2.
\end{eqnarray*}
\end{proposition}
{\it Proof}. By Gaffney's integration by parts formula, we have
\begin{eqnarray*}
\|dd^*(a+\square)^{-1}\omega\|_2^2&=&\int_M \langle dd^*(a+\square)^{-1}\omega, dd^*(a+\square)^{-1}\omega\rangle dv\\
&=&\int_M\langle (a+\square)^{-1}\omega, dd^*dd^*(a+\square)^{-1}\omega\rangle dv.
\end{eqnarray*}
Similarly, we can prove
\begin{eqnarray*}
\|d^*d(a+\square)^{-1}\omega\|_2^2=\int_M\langle (a+\square)^{-1}\omega, d^*dd^*d(a+\square)^{-1}\omega\rangle dv.
\end{eqnarray*}
Using the fact that  $dd^*dd^*+d^*dd^*d=\square^2$, we get
\begin{eqnarray*}
\|dd^*(a+\square)^{-1}\omega\|_2^2+\|d^*d(a+\square)^{-1}\omega\|_2^2=\int_M \langle(a+\square)^{-1}\omega, \square^2(a+\square)^{-1}\omega\rangle dv.
\end{eqnarray*}
This proves the identity $(\ref{ID1})$. Again, integration by parts yields
\begin{eqnarray*}
\|(a+\square)\omega\|_2^2&=&\|\square\omega\|_2^2+2a\langle \langle\omega, \square\omega\rangle\rangle+a^2\|\omega\|_2^2\\
&=&\|\square\omega\|_2^2+2a\|d\omega\|_2^2+2a\|d^*\omega\|_2^2+a^2\|\omega\|_2^2\\
&\geq& \|\square\omega\|_2^2,
\end{eqnarray*}
which implies that
\begin{eqnarray}
\|\square(a+\square)^{-1}\omega\|_2\leq \|\omega\|_2.\label{ID2}
\end{eqnarray}
Combining $(\ref{ID1})$ with $(\ref{ID2})$, we obtain
\begin{eqnarray*}
\|dd^*(a+\square)^{-1}\omega\|_2^2+\|d^*d(a+\square)^{-1}\omega\|_2^2\leq \|\omega\|_2^2.
\end{eqnarray*}
This finishes the proof of Proposition \ref{BA2}. \hfill $\square$

\medskip

We now state the $L^p$-norm estimates of the Beurling-Ahlfors transforms on complete Riemannian manifolds. The following result is the restatement of Theorem 1.2 and Theorem 5.1 in \cite{Li2011}. Here, as in \cite{Li2011}, $\|\cdot\|_{\rm op}$ denotes the operator norm.

\begin{theorem}\label{BA}
Suppose that there exists a constant $a\geq 0$ such that
$$W_k\geq -a.$$
Then, there exists a universal constant $C>0$ such that for all $1<p<\infty$, and for all $\omega\in C_0^\infty(\Lambda^kT^*M)$,
\begin{eqnarray*}
\|S_{A_i}\omega\|_p\leq C(p^*-1)^{3/2}\|A_i\|_{\rm op}\|\omega\|_p,
\end{eqnarray*}
and
\begin{eqnarray*}
\|S_{B}\omega\|_p\leq C(p^*-1)^{3/2}\|B\|_{\rm op}\|\omega\|_p.
\end{eqnarray*}
In particular, in the case where $W_k\equiv -a$, we have
\begin{eqnarray*}
\|S_{A_i}\omega\|_p\leq 2(p^*-1)\|A_i\|_{\rm op}\|\omega\|_p,
\end{eqnarray*}
and
\begin{eqnarray*}
\|S_{B}\omega\|_p\leq 2(p^*-1)\|B\|_{\rm op}\|\omega\|_p.
\end{eqnarray*}
\end{theorem}
{\it Proof}. By Proposition \ref{BA2}, we need only to study the case $p\neq 2$. For simplicity, we only consider the case $W_k\geq 0$. The general case $W_k\geq -a$ can be similarly proved. Let
\begin{eqnarray*}
Z_{t}^{i}=M_t\int_0^t M_s^{-1}A_i\nabla \omega(X_s, t-s)dX_s, \ \ \ i=1, 2.
\end{eqnarray*}
By Theorem 2.6 due to Ba\~nuelos and Baudoin in \cite{BB}, for all $p\in (1, \infty)$, we have
\begin{eqnarray*}
\|Z_{T}^{i}\|_p\leq 3\sqrt{p(2p-1)}
\left\|\left(\int_0^T|A_i\nabla\omega(X_t, T-t)|^2dt\right)^{1\over
2}\right\|_p.
\end{eqnarray*}
Obviously, we have

\begin{eqnarray*}
\left\|\left(\int_0^T|A_i\nabla\omega(X_t, T-t)|^2dt\right)^{1\over
2}\right\|_p\leq \|A_i\|_{\rm
op}\left\|\left(\int_0^T|\nabla\omega(X_t, T-t)|^2dt\right)^{1\over
2}\right\|_p.
\end{eqnarray*}
By the same argument as used in the proofs of Proposition 6.2 and Proposition 6.3 in \cite{Li2010}, for all $1<p<\infty$, we can prove that
\begin{eqnarray*}
\left\|\left(\int_0^T|\nabla\omega(X_t, T-t)|^2dt\right)^{1\over
2}\right\|_p\leq B_p\|\omega\|_p,
\end{eqnarray*}
where $B_p=(2p)^{1/2}(p-1)^{-3/2}$ for $p\in (1, 2)$, $B_p=1$ for $p=2$, and
$B_p={p\over \sqrt{2(p-2)}}$ if $p>2$.
Hence, for $1<p<2$,
\begin{eqnarray*}
\|S_{A_i}^T\omega\|_p&\leq &3\sqrt{p(2p-1)} \|A_i\|_{\rm op} {(2p)^{1/2}\over (p-1)^{3/2}}\|\omega\|_p\\
&\leq& 6\sqrt{6}(p-1)^{-3/2}\|A_i\|_{\rm op}\|\omega\|_p,
\end{eqnarray*}
and for $p>2$,
\begin{eqnarray*}
\|S_{A_i}^T\omega\|_{p}\leq 3(p-1)^{3/2}(1+O((p-1)^{-1}))\|A_i\|_{\rm op}\|\omega\|_p.
\end{eqnarray*}
Indeed, by duality argument as used in  \cite{Li2011}, for all $p>2$,
we have
\begin{eqnarray*}
\|S_{A_i}^T\|_{p, p}=\|S_{A_i}^T\|_{q, q},
\end{eqnarray*}
which yields
for $p>2$,
\begin{eqnarray*}
\|S_{A_i}^T\omega\|_{p}\leq 6\sqrt{6}(p-1)^{3/2}\|A_i\|_{\rm op}\|\omega\|_p.
\end{eqnarray*}
In summary, for all $1<p<\infty$, we have proved that
\begin{eqnarray*}
\|S_{A_i}^T\omega\|_{p}\leq 6\sqrt{6}(p^*-1)^{3/2}\|A_i\|_{\rm op}\|\omega\|_p.
\end{eqnarray*}
Similarly, for all $1<p<\infty$, we can prove
\begin{eqnarray*}
\|S_{B}^T\omega\|_{p}\leq 6\sqrt{6}(p^*-1)^{3/2}\|B\|_{\rm op}\|\omega\|_p.
\end{eqnarray*}

In the particular case where $W_k\equiv -a$, we have
\begin{eqnarray*}
S_{A_i}^T\omega(x)=E\left[\left. U_T\int_0^T U_t^{-1}A_i\nabla
\omega_a(X_t, T-t)dX_t\right|X_T=x\right], \ \ \ i=1, 2, \ {\rm
a.s.} x\in M.
\end{eqnarray*}
The $L^p$-contractiveness of the conditional expectation yields

\begin{eqnarray*}
\|S_{A_i}^T\omega\|&=&\left\|U_T\int_0^T U_t^{-1}A_i\nabla
\omega_a(X_t, T-t)dX_t\right\|_p\\
&=&\left\|\int_0^T U_t^{-1}A_i\nabla \omega_a(X_t, T-t)dX_t\right\|_p.
\end{eqnarray*}
Using the Burkholder sharp $L^p$-inequality for martingale
transforms, for all $p>1$, we deduce that
\begin{eqnarray}
\|S_{A_i}^T\omega\|\leq (p^*-1)\sup\limits_{0\leq t\leq
T}\|U_t^{-1}A_iU_t\|_{\rm op}\left\|\int_0^T U_t^{-1}\nabla
\omega_a(X_t, T-t)dX_t\right\|_p.\label{SA-8}
\end{eqnarray}
By It\^o's formula, we can prove that (see Eq. $(49)$ in \cite{Li2011})
\begin{eqnarray}
\omega(X_T)-U_T\omega_a(X_0, T)=U_T\int_0^T U_t^{-1}\nabla\omega_a(X_t,
T-t)dX_t.\label{SA-9}
\end{eqnarray}
Substituting $(\ref
{SA-9})$ into $(\ref{SA-8})$, we have
\begin{eqnarray*}
\|S_{A_i}^T\omega\|_p\leq (p^*-1)\|A_i\|_{\rm op}
\|\omega(X_T)-U_T\omega_a(X_0, T)\|_p.
\end{eqnarray*}
Using the argument in \cite{Li2011}, we obtain
\begin{eqnarray*}
\|S_A^T\omega\|_p\leq (p^*-1)\left(1+e^{-2\min\{{1\over p},
1-{1\over p}\}aT}\right)\|\omega\|_p.
\end{eqnarray*}
Hence
\begin{eqnarray*}
\|dd^*(a+\square)^{-1}\omega\|_p\leq 2\lim\limits_{T\rightarrow
\infty}\|S_{A_1}^T\omega\|_p\leq 2(p^*-1)\|A_1\|_{\rm
op}\|\omega\|_p,
\end{eqnarray*}
\begin{eqnarray*}
\|d^*d(a+\square)^{-1}\omega\|_p\leq 2\lim\limits_{T\rightarrow
\infty}\|S_{A_2}^T\omega\|_p\leq 2(p^*-1)\|A_2\|_{\rm
op}\|\omega\|_p,
\end{eqnarray*}
and
\begin{eqnarray*}
\|(dd^*-d^*d)(a+\square)^{-1}\omega\|_p\leq
2\lim\limits_{T\rightarrow \infty}\|S_{B}^T\omega\|_p\leq
2(p^*-1)\|B\|_{\rm op}\|\omega\|_p.
\end{eqnarray*}
The proof of Theorem \ref{BA} is completed. \hfill $\square$

\begin{remark}\label{rrr} {\rm The above proof corrects a gap contained in \cite{Li2011}.  In summary, the $L^p$-norm estimates in  Theorem \ref{BA} indicates that the results in Theorem 1.2, Theorem 1.3 and Theorem 1.4, Theorem 5.1 and Corollary 5.2 obtained in \cite{Li2011} remain valid. As a consequence, the main theorems proved in \cite{Li2011} remain valid. In particular, see Theorem 1.3 in \cite{Li2011},  on  complete and stochastically complete Riemannian manifolds non-negative Weitzenb\"ock curvature operator $W_k\geq 0$, where $1\leq k\leq n={\rm dim}M$, the Weak $L^p$-Hodge decomposition theorem holds for $k$-forms, the De Rham projection ${\rm P}_1=dd^*\square^{-1}$, the Leray projection ${\rm P}_2=d^*d\square^{-1}$ and  the Beurling-Ahlfors transform $B_k=(d^*d-dd^*)\square^{-1}$ on $k$-form is bounded in $L^p$ for all $1<p<\infty$.}

\end{remark}

\section{Time reversal martingale transformation representation formula for the Riesz transfroms}

In this section, we prove a time reversal martingale transformation representation formula
for the Riesz transforms on complete Riemannian
manifolds.

\medskip
First, we prove the following time reversal martingale transformation representation formula for one forms.

\begin{theorem}\label{theo1-b} Let $\widehat{X}_t=X_{\tau-t}$, and $\widehat{B}_t=B_{\tau-t}$, $t\in [0,
\tau]$. Let $\widehat{M}_t$ be the solution to the covariant SDE
\begin{eqnarray*}
{\nabla\over \partial t}\widehat
{M}_t&=&-\widehat{M}_t(Ric+\nabla^2\phi)(\widehat{X}_t),\\
\widehat{M}_0&=&{\rm Id}_{T_{\widehat X_0}M}.
\end{eqnarray*}
For any $\omega\in C_0^\infty(\Lambda^1 T^*M)$, let $\omega_a(x,
y)=e^{-y\sqrt{a+\square_\phi}}\omega(x)$, $\forall x\in M, y\geq 0$.
Then, for a.s. $x\in M$,
\begin{eqnarray*}
{1\over 2}\omega(x)&=&\lim\limits_{y\rightarrow +\infty}E_y\left[
\left. \widehat Z_\tau \right|\widehat{X}_0=x\right],\label{F2-b}
\end{eqnarray*}
where
\begin{eqnarray*}
\widehat Z_\tau=\int_0^\tau
e^{-at}\widehat{M_t}\partial_y\omega_a(\widehat{X}_t,
\widehat{B}_t)d\widehat{B}_t-\int_0^\tau e^{-at}
\widehat{M}_{t}\partial_y^2 \omega(\widehat{X}_t, \widehat{B}_t)dt.
\end{eqnarray*}
\end{theorem}
{\it Proof}. By Theorem \ref{theo1}, we have
\begin{eqnarray*}
{1\over 2}\omega(x)&=&\lim\limits_{y\rightarrow +\infty}E_y\left[
\left. Z_\tau \right|X_\tau=x\right],
\end{eqnarray*}
where
$$Z_\tau=e^{-a\tau}M_\tau\int_0^\tau e^{as}M_s^{-1}\nabla_y \omega_a(X_s,
B_s)dB_s.$$ Taking $0=s_0<s_1<\ldots<s_n<s_{n+1}=\tau$ be a partition
of $[0, \tau]$, then
\begin{eqnarray*}
Z_{\tau, n}:= e^{-a\tau} M_\tau\sum\limits_{i=1}^N e^{as_i}M_{s_i}^{-1}\nabla_y \omega(X_{s_i}, B_{s_i})(B_{s_{i+1}}-B_{s_i})
\end{eqnarray*}
converges in $L^2$ and in probability to $Z_\tau$. We can rewrite $Z_{\tau, n}$ as follows
\begin{eqnarray*}
Z_{\tau, n}=\sum\limits_{i=1}^N e^{-a(\tau-s_i)}M_\tau
M_{s_i}^{-1}\nabla_y \omega(X_{s_i}, B_{s_i})(B_{s_{i+1}}-B_{s_i}).
\end{eqnarray*}

Note that
\begin{eqnarray*}
\partial_s\widehat{M}_{\tau-s}&=&\widehat{M}_{\tau-s} Ric(L)(\widehat{X}_{\tau-s})\\
&=&\widehat{M}_{\tau-s} Ric(L)(X_{s}),
\end{eqnarray*}
and
\begin{eqnarray*}
\partial_s (M_{\tau}M_{s}^{-1})&=&-M_{\tau}M_s^{-1}\partial_s M_s M_s^{-1}\\
&=&M_{\tau}M_{s}^{-1} Ric(L)(X_s)M_s M_s^{-1}\\
&=&(M_{\tau}M_s^{-1}) Ric(L)(X_s).
\end{eqnarray*}
By the uniqueness of the solution to ODE, as
$\left.\widehat{M}_{\tau-s}\right|_{s=\tau}=\left.M_\tau
M_{s}^{-1}\right|_{s=\tau}={\rm Id}_{T_{\widehat{X}_0}M}$,  we have
\begin{eqnarray*}
M_\tau M_s^{-1}=\widehat{M}_{\tau-s}.
\end{eqnarray*}
Therefore
\begin{eqnarray*}
Z_{\tau, n}=\sum\limits_{i=1}^N e^{-a(\tau-s_i)}\widehat{M}_{\tau-s_i}\nabla_y \omega(\widehat{X}_{\tau-s_i}, \widehat{B}_{\tau-s_i})(\widehat{B}_{\tau-s_{i+1}}-\widehat{B}_{\tau-s_i})\\
\end{eqnarray*}
Let $t_i=\tau-s_i$. Then $\tau=t_0>t_1>\ldots>t_n>t_{n+1}=0$, and
\begin{eqnarray*}
Z_{\tau, n}=\sum\limits_{i=1}^N e^{-at_i}\widehat{M}_{t_i}\nabla_y \omega(\widehat{X}_{t_i}, \widehat{B}_{t_i})(\widehat{B}_{t_{i+1}}-\widehat{B}_{t_{i+1}}).
\end{eqnarray*}
By Taylor's formula, we have
\begin{eqnarray*}
\omega(\widehat{X}_{t_i}, \widehat{B}_{t_i})=
\omega(\widehat{X}_{t_{i+1}}, \widehat{B}_{t_{i+1}}) -\nabla_y
\omega(\widehat{X}_{t_{i+1}},
\widehat{B}_{t_{i+1}})(\widehat{B}_{t_{i+1}}-\widehat{B}_{t_{i}})+O\left((\widehat{B}_{t_{i+1}}-\widehat{B}_{t_{i}})^2\right).
\end{eqnarray*}
Hence
\begin{eqnarray*}
Z_{\tau, n}&=&\sum\limits_{i=1}^N e^{-at_i}\widehat{M}_{t_i}\nabla_y \omega(\widehat{X}_{t_{i+1}}, \widehat{B}_{t_{i+1}})(\widehat{B}_{t_{i+1}}-\widehat{B}_{t_i})\\
& &\ \ \ -\sum\limits_{i=1}^N e^{-at_i}\widehat{M}_{t_i}\nabla_y^2
\omega(\widehat{X}_{t_{i+1}},
\widehat{B}_{t_{i+1}})(\widehat{B}_{t_{i+1}}-\widehat{B}_{t_i})^2+O\left((\widehat{B}_{t_{i+1}}-\widehat{B}_{t_{i}})^3\right).
\end{eqnarray*}
which converges in $L^2$ and in
probability to the following limit
\begin{eqnarray*}
 \widehat Z_{\tau}=\int_0^\tau e^{-at}\widehat{M}_{t}\nabla_y \omega(\widehat{X}_{t}, \widehat{B}_t)d\widehat{B}_t-\int_0^\tau e^{-at}
 \widehat{M}_{t}\nabla_y^2 \omega(\widehat{X}_t, \widehat{B}_t)dt.
\end{eqnarray*}
The proof of Theorem \ref{theo1-b} is completed. \hfill $\square$

By Theorem \ref{theo1-b}, we can prove the following time reversal martingale
transformation representation formula for the Riesz transforms on
complete Riemannian manifolds.

\begin{theorem} Let $R_a(L)= \nabla(a-L)^{-1/2}$. Then, for $f\in C_0^\infty(M)$, we have
\begin{eqnarray*}
R_a(L)f(x)=-2\lim\limits_{y\rightarrow +\infty}E_y\left[
\left. \widehat Z_\tau \right|\widehat{X}_0=x\right]\label{R1-b},
\end{eqnarray*}
where
\begin{eqnarray*}
\widehat Z_\tau=\int_0^\tau e^{-as}\widehat M_s dQ_af(\widehat{X}_s,
\widehat{B}_s)d\widehat{B}_s-\int_0^\tau e^{-as}\widehat M_s\partial_ydQ_af(\widehat{X}_s, \widehat{B}_s)ds.
\end{eqnarray*}
\end{theorem}

\begin{remark} {\rm
As noticed in \cite{Li2008}, there exists a standard one dimensional Brownian motion $\beta_t$ such that
\begin{eqnarray*}
d\widehat{B}_t=d\beta_t+{dt\over \widehat{B}_t}, \ \ \ t\in (0, \tau].
\end{eqnarray*}}
\end{remark}

\section{Time reversal martingale transforms representation formula for the Beurling-Ahlfors transforms}

Similarly to the proof of Theorem \ref{theo1-b},  we prove a time reversal martingale transformation representation formula
for the Beurling-Ahlfors transforms on complete Riemannian manifolds.

\begin{theorem}  Let $\widehat X_t=X_{T-t}$, $t\in [0, T]$. Let $\widehat M_t$ be the solution to the covariant equation
$${\nabla \widehat M_t\over \partial t}=-\widehat M_t W_k(\widehat X_t), \ \ \widehat M_0={\rm Id}_{\Lambda^kT^*_{\widehat X_0}M}.$$
Then, for any $\omega\in C_0^\infty(\Lambda^1 T^*M)$, the Beurling-Ahlfors transform
$$
S_B\omega:=(d^*d-dd^*)(a+\square)^{-1}\omega
$$
has the following time reversal martingale transform representation: for a.s. $x\in M$,
\begin{eqnarray*}
S_B\omega(x)=2\lim\limits_{T\rightarrow \infty}E\left[\left. \widehat Z_T\right|\widehat X_0=x\right],
\end{eqnarray*}
where
\begin{eqnarray*}
\widehat Z_T=\int_0^T e^{-as }\widehat M_s B \nabla \omega_a(\widehat X_s, s)d\widehat X_s-\int_0^T e^{-as} \widehat M_s B{\rm Tr}\nabla^2\omega_a(\widehat X_s, s)ds.
\end{eqnarray*}

\end{theorem}

To end this paper, let us mention that, in a forthcoming paper \cite{Li2013},  we will prove a martingale transform representation formula for the Riesz transforms associated with the Dirac operator acting on Hermitian vector bundles over complete Riemannian manifolds and for the Riesz transforms associated with the $\bar\partial$-operator acting on holomorphic Hermitian vector bundles over complete K\"ahler manifolds. By the same argument as used in this paper, we can prove some explicit dimension free $L^p$-norm estimates of these Riesz transforms on complete Riemannian or K\"ahler manifolds with suitable curvature conditions. See also \cite{Li2010b}.

\medskip

\noindent{\bf Acknowledgement}.   \ \ I would like to thank R.
Ba\~nuelos and F. Baudoin for their interests on my previous works
and for pointing out the gap contained in \cite{Li2008, Li2011}
which has been addressed in this paper.

\end{document}